\documentclass[preprint,12pt]{elsarticle}

\usepackage{amsmath,amsfonts,latexsym,subfigure,amsthm,amssymb,wasysym}
\usepackage{graphicx}
\usepackage{tikz}
\usepackage{ifthen}
\usepackage[usenames,dvipsnames]{pstricks}
\usepackage{epsfig}
\usepackage{pst-grad} 
\usepackage{pst-plot} 
\usepackage{enumitem}
\usepackage{algorithm}
\usepackage{algorithmic}
\usepackage{wrapfig}
\usepackage{color}

\def\calN{{\mathcal N}}
\def\calD{{\mathcal D}}

\def\calC{{\mathcal C}}

\def\calV {{\mathcal V}}
\def\calY {{\mathcal Y}}
\def\calJ {{\mathcal J}}
\def\calB {{\mathcal B}}

\def\calE {{\mathcal E}}
\def\calR {{\mathcal R}}

\def\mubar {{\bar{\boldsymbol \mu}}}
\def\vmu {{\boldsymbol \mu}}
\newcommand{\x}{\boldsymbol{x}}
\def\betabar {{\bar\beta}}

\usepackage{bbm}

\theoremstyle{remark}

\journal{Applied Numerical Mathematics}

\begin{document}

\begin{frontmatter}



\title{A Certified Natural-Norm Successive Constraint Method for Parametric Inf-Sup Lower Bounds}


\author[UMassD]{Yanlai Chen \fnref{funding}}
\ead{{yanlai.chen@umassd.edu}}
\ead[url]{www.faculty.umassd.edu/yanlai.chen}
\address[UMassD]{Department of Mathematics, University of Massachusetts Dartmouth, 285 Old Westport Road, North Dartmouth, MA 02747, USA.}

\fntext[funding]{This research was partially supported by National Science Foundation grant DMS-1216928.}

\begin{abstract}
We present a certified version of the Natural-Norm Successive Constraint Method (cNNSCM) 
for fast and accurate Inf-Sup lower bound evaluation of parametric operators. Successive Constraint Methods (SCM) are essential tools 
for the construction of a lower bound for the inf-sup stability constants which are required in  {\it a posteriori} error analysis 
of reduced basis approximations. 
They utilize a Linear Program (LP) relaxation scheme 
incorporating continuity and stability constraints. 
The natural-norm approach {\em linearizes} inf-sup constant as a function of the parameter. 
The Natural-Norm Successive Constraint Method (NNSCM) combines these two aspects. It uses a greedy algorithm to select 
SCM control points which adaptively construct an optimal decomposition of the parameter domain, and then 
apply the SCM on each domain. 

Unfortunately, the NNSCM produces no guarantee for the quality of the lower bound. The new cNNSCM provides 
an upper bound in addition to the lower bound and let the user control the gap, thus the quality of the lower bound.
The efficacy and accuracy of the new method is validated by numerical experiments.
\end{abstract}


\end{frontmatter}



\section{Introduction}
\label{sec:intro}

For affinely parametrized partial differential equations, the certified reduced basis method (RBM) 
\cite{NguyenVeroyPatera2005,Prudhomme_Rovas_Veroy_Maday_Patera_Turinici,Rozza_Huynh_Patera, HaasdonkOhlberger} 
utilizes an Offline-Online computational decomposition strategy to produce surrogate solution (of dimension $N$) in a 
time that is of orders of magnitude shorter than what is needed by the underlying numerical solver of dimension ${\mathcal N} \gg N$ (called 
{\em truth} solver hereafter). 
The RBM relies on a projection onto a low dimensional space spanned by truth approximations at 
an optimally sampled set of parameter values \cite{
Almroth_Stern_Brogan,Fink_Rheinboldt_1,Noor_Peters,Porsching, Maday}. 
This low-dimensional manifold is generated by a greedy algorithm making use of a rigorous {\it a posteriori} error bounds
for the field variable and associated functional outputs of interest which also guarantees the fidelity 
of the surrogate solution in approximating the truth approximation. 
The high efficiency and accuracy of RBM 
render it an ideal candidate for practical methods in the real-time and many-query contexts. 

This crucial {\it a posteriori} error bound is residual-based and requires an estimate (lower bound) 
for the stability factor of the discrete partial differential operator, that is the coercivity 
or inf-sup constant. 
In the RBM context, given any parameter value this stability factor must be estimated efficiently. 
So it should also admit an Offline-Online structure for which the Online expense is independent of $\mathcal N$. 
Moreover, the optimality of the low-dimensional RB manifold is dependent on the quality 
of this estimate as a parameter-dependent function, so the lower bound should not be too pessimistic. 
There are several approaches in the literature. 
A Successive Constraint Method (SCM) is proposed in 
\cite{HuynhSCM} and subsequently improved in \cite{CHMR-Cras,CHMR-M2an,Vallaghe,Zhang2011}.
It is a framework incorporating both continuity {\em and} stability information whose Online component 
is the resolution of a small-size Liner Programming (LP) problem. Hence, this procedure is rather efficient.
However, the classical inf-sup formulation 
has couple of undesirable attributes -- a $Q^2$-term affine parameter expansion (resulting from a squaring of the operator), 
and loss of (even local) concavity.
On the other hand, a ``natural-norm" method is proposed in \cite{Deparis08,SenNatNorm}. Its linearized-in-parameter 
inf-sup formulation has several
desirable approximation properties - a $Q$-term affine parameter expansion, and first order (in parameter) concavity;
however, the lower bound procedure is rather crude - a framework which incorporates only continuity information.
A natural-norm SCM approach is proposed in \cite{HKCHP} combining the ``linearized" inf-sup statement 
with the SCM lower bound procedure. 
The former (natural-norm) provides 
a smaller optimization problem which enjoys intrinsic lower bound properties. 
The latter (SCM) provides a systematic optimization framework: a Linear Program  
relaxation which readily incorporates effective stability constraints.
The natural-norm SCM  performs very well in
particular in the Offline stage:  it is typically an order of magnitude less costly than either
the natural-norm or ``classical" SCM approaches alone. However, unlike the classical SCM, it provides no 
upper-bound thus no control of the quality of the lower bound. This often results in extremely pessimistic estimate. 

In this paper, we propose a certified version of the NNSCM (cNNSCM). Without significantly degrading the efficiency, it 
provides an upper-bound and thus a mechanism for the user to easily control the quality of the lower bound. 
As a result, the lower bound of the new cNNSCM may be orders of magnitude more accurate than the original NNSCM. 
The method is tested on two elliptic partial differential equations. In what follows, we use the same notation as 
in \cite{HKCHP} and 
denote the classical SCM method \cite{HuynhSCM,CHMR-Cras,CHMR-M2an} 
as $\rm SCM^2$ in order to differentiate it from the new 
natural-norm type of approaches. Here, the (squared) superscript suggests the undesired $Q^2$-term 
affine parameter expansion in the classical method.

This paper is organized as follows. In Section \ref{sec:background}, we review the background materials including 
the RBM, its {\em A Posteriori} error estimation and the involved stability constant. Section \ref{sec:lower-bound} 
describes the natural-norm SCM. The new certified NNSCM is proposed in Section \ref{sec:cnnscm}. Numerical 
validations are presented in Section \ref{sec:numerical}, and finally some concluding remarks are offered in Section \ref{sec:conclusion}.

\section{Background}

\label{sec:background}

For the completeness of this paper and to put the concerned method into context, we introduce the necessary background 
materials in this section. 
To that end, this section covers the truth solver and the related stability constants, the reduced basis method, 
and the {\it A Posteriori} error estimate needed therein.

\label{sec:formulation}

\subsection{Notations}
We use $\Omega \subset {\mathbb R}^n$ ($ n = {2 \mbox{ or } 3}$) to denote 
a bounded physical domain with boundary $\partial \Omega$. 
We introduce a closed parameter domain $\calD \in {\mathbb R}^P$, 
a point ($P$-tuple) in which is denoted ${\boldsymbol \mu} = (\mu_1,\ldots,\mu_P)$. A set of $N$ parameter values will be differentiated by 
superscripts $\{{\boldsymbol{\mu}}^i\}_{i=1}^N$.
Let us then define the Hilbert space
$X$  equipped with inner product $(\cdot,\cdot)_X$
and induced norm $\| \cdot \|_X$. Here $(H^1_0(\Omega))^\calV \subset X \subset (H^1(\Omega))^\calV$
($\calV = 1$ for a scalar field and $\calV > 1$ for a vector field) 
\cite{QuarteroniVallie2008book, Adams1975Book}. 
Finally, we introduce a parametrized bilinear form  and two linear forms. 
$a(\cdot,\cdot;{\boldsymbol \mu})$: $X \times X \rightarrow \mathbb{R}$ is such that 
\begin{itemize}
\item it is
inf-sup stable and continuous over $X$: $\beta({\boldsymbol \mu}) > 0$ and $\gamma({\boldsymbol \mu})$ is
finite $\forall {\boldsymbol \mu} \in \calD$, where
$$ \beta({\boldsymbol \mu}) = \inf_{w\in X} \sup_{v \in X} \frac{a(w,v;{\boldsymbol \mu})}{\|w\|_X\,\|v\|_X} , \mbox{ and } \gamma({\boldsymbol \mu}) = \sup_{w \in X} \sup_{v \in X} \frac{a(w,v;{\boldsymbol \mu})}{\|w\|_X\,\|v\|_X} ;$$ 
\item $a$ is ``affine" in the parameter: 
$a(w,v;{\boldsymbol \mu}) = {\displaystyle \sum_{q = 1}^Q} \Theta_q({\boldsymbol \mu}) a_q(w,v)$.
\end{itemize}
We emphasize that it can be
approximated by affine (bi)linear forms when it is nonaffine
\cite{Barrault_Nguyen_Maday_Patera,Grepl_Maday_Nguyen_Patera}.
Finally, we introduce two linear bounded functionals $f(\cdot; {\boldsymbol \mu}): X \rightarrow {\mathbb R}$ and $\ell(\cdot;{\boldsymbol \mu}): X \rightarrow \mathbb{R}$ 
that are also affine in the parameter.
The following continuous problem is then well-defined. 

{\bf ($P_C$)} Given ${\boldsymbol \mu} \in \calD$,  find $u({\boldsymbol \mu}) \in X$ such that
$ a(u({\boldsymbol \mu}),v;{\boldsymbol \mu}) = f(v, \vmu), \forall v \in X$.

For many applications, we concern a scalar quantify of interest as $s({\boldsymbol \mu}) = \ell(u({\boldsymbol \mu}), \vmu)$.  
To discretize this problem, we consider for an example a finite element approximation space (of dimension $\calN$) $X^\calN \subset X$. 
Suppose that the discretized bilinear form 
remains inf-sup stable (and continuous) over $X^\calN$ with constants $\beta^\calN({\boldsymbol \mu}) > 0$ and
$\gamma^\calN({\boldsymbol \mu})$ being finite $\forall {\boldsymbol \mu} \in \calD$, where
$$ \beta^\calN({\boldsymbol \mu}) = \inf_{w\in X^\calN} \sup_{v \in X^\calN} \frac{a^\calN(w,v;{\boldsymbol \mu})}{\|w\|_{X^\calN}\,\|v\|_{X^\calN}}  \mbox{ and }
\gamma^\calN({\boldsymbol \mu}) = \sup_{w \in X^\calN} \sup_{v \in X^\calN} \frac{a^\calN(w,v;{\boldsymbol \mu})}{\|w\|_{X^\calN}\,\|v\|_{X^\calN}} .$$
We now introduce our {\em truth} discretization: 

{\bf ($P_D$)} Given ${\boldsymbol \mu} \in \calD$,  find $u^\calN({\boldsymbol \mu}) \in X^\calN$ such that $ a(u^\calN({\boldsymbol \mu}),v;{\boldsymbol \mu}) = f(v, \vmu), \forall v \in X^\calN$.

This discretization is called {\em truth} and its solution {\em truth approximation} because 
our reduced basis approximation is built upon, and its error measured  with 
respect to this finite element solution. 
The (numerical) output is evaluated accordingly $s^\calN({\boldsymbol \mu}) = \ell(u^\calN({\boldsymbol \mu}), \vmu)$.  
Since we are essentially abandoning {\bf ($P_C$)}, for brevity of exposition we may omit the $\calN$ when there is no confusion. 

We end this section by re-writing the inf-sup constant $\beta({\boldsymbol \mu})$.
To that end, we first define the supremizing operator $T^{\boldsymbol \mu}: X \rightarrow X$ such that $(T^{\boldsymbol \mu} w, v)_X = a(w,v;{\boldsymbol \mu}), \forall v \in X$. Clearly, we have 
$$ T^{\boldsymbol \mu} w = \arg \sup_{v \in X} \frac{a(w,v;{\boldsymbol \mu})}{\|v\|_X}, \mbox{ and that furthermore }
\beta({\boldsymbol \mu}) = \inf_{w \in X}   \frac{ \| T^{\boldsymbol \mu} w\|_X }{ \| w \|_X }. $$
Recalling the affine assumption allows us to decompose $T^{{\boldsymbol \mu}}$ as
$$ T^{\boldsymbol \mu} w = \sum_{q = 1}^Q \Theta_q({\boldsymbol \mu}) T_q w$$
where $(T_q w, v)_X = a_q(w,v), \forall v \in X$, $1 \le q \le Q$.

\subsection{Reduced Basis Method and the \emph{A Posteriori} Error Estimators}
\label{subsec:error-estimators}

The fundamental observation utilized by RBM 
is that $u^{\mathcal N}({\boldsymbol \mu})$ residing on ${\mathcal M} = \{u^{\mathcal N}({\boldsymbol \mu}),\,{\boldsymbol \mu} \in {\mathcal
D}\}$ can typically be well approximated by a finite-dimensional space. 
The RBM idea is then to propose an approximation  of ${\mathcal M}$ by 
$$W^N={\rm
span}\{u^{\mathcal N}({\boldsymbol \mu}^1),\,\dots,\,u^{\mathcal N}({\boldsymbol \mu}^N)\}$$ where, $u^{\mathcal
N}({\boldsymbol \mu}^1),\,\dots,\,u^{\mathcal N}({\boldsymbol \mu}^N)$ are $N$ $(\ll\mathcal N)$ truth approximations corresponding to
the parameters $\{{\boldsymbol \mu}^1,\dots,{\boldsymbol \mu}^N\}$ selected according to a judicious sampling strategy \cite{Maday}. For
a given ${\boldsymbol \mu}$, we now 
solve in $W^N$ for the reduced solution $u^{N}({\boldsymbol \mu})$. 

{\bf ($P_R$)} Given ${\boldsymbol \mu} \in \calD$,  find $u^N({\boldsymbol \mu}) \in W^N$ such that $ a(u^N({\boldsymbol \mu}),v;{\boldsymbol \mu}) = f(v), \forall v \in W^N$.

The online computation is $\mathcal N$-independent, thanks to the
assumption that the (bi)linear forms are affine. Hence, the online part is very efficient. In
order to be able to ``optimally'' find the $N$ parameters and to assure the fidelity of the reduced basis
solution $u^{N}({\boldsymbol \mu})$ to approximate the truth solution $u^{\mathcal N}({\boldsymbol \mu})$, we need 
the {\em a posteriori} error estimator $\Delta_N({\boldsymbol \mu})$ 
\cite{Machielis_Maday_Oliveira_Patera_Rovas,Maday_Patera_Rovas,Prudhomme_Rovas_Veroy_Maday_Patera_Turinici,
Rozza_Huynh_Patera,SenNatNorm} 
that involves the residual $$r(v,{\boldsymbol \mu}) \equiv f^\calN(v;{\boldsymbol \mu}) - a^\calN (u^N({\boldsymbol \mu}),v;{\boldsymbol \mu})$$ and the inf-sup stability constant $\beta^\calN({\boldsymbol \mu})$.
With this estimator, we can describe briefly the classical greedy algorithm used to find the $N$ parameters
${\boldsymbol \mu}_1, \dots, {\boldsymbol \mu}_N$ and the space $W^N$: 
We first randomly
select one parameter value and compute the associated truth approximation. Next, we scan the entire discrete parameter space and for each parameter in this space compute its RB approximation $u^{N=1}$ and the error estimator
$\Delta_1({\boldsymbol \mu})$. The next parameter value we select, ${\boldsymbol \mu}_2$, is the one corresponding to the largest error
estimator. We then compute the truth approximation and thus have a new basis set
consisting 
of two elements.
This process is repeated until the maximum of the error estimators is sufficiently small.

We end by providing the missing component - how the inf-sup lower bound will serve within the error 
estimators. 
The reduced basis field error and output error (relative to the truth discretization) satisfies \cite{NguyenVeroyPatera2005,SenNatNorm}
\begin{itemize}
\item $ | u^{\calN}({\boldsymbol \mu}) - u^N({\boldsymbol \mu}) | \le \Delta_N({\boldsymbol \mu}),\mbox{ where } \Delta_N({\boldsymbol \mu}) \equiv \frac { \| r(\cdot;{\boldsymbol \mu}) \|_{X'} } {\beta^{\rm LB}({\boldsymbol \mu})},$
\item $ | s^{\calN}({\boldsymbol \mu}) - s_N({\boldsymbol \mu}) | \le \Delta^s_N({\boldsymbol \mu}),\mbox{ where } \Delta^s_N({\boldsymbol \mu}) \equiv \frac { \|\ell (\cdot) \|_{X'} \| r(\cdot;{\boldsymbol \mu}) \|_{X'} } {\beta^{\rm LB}({\boldsymbol \mu})}.$
\end{itemize}
Here, $\| \cdot \|_{X'}$ refers to the dual norm with respect to $X^{\calN}$ and $
\beta^{\rm LB}$ is a lower bound of $\beta^\calN({\boldsymbol \mu})$. 
The later implies that the quality of the inf-sup lower bound affects the quality of the error bound which, in turn, 
determines the optimality of the RB space $W^N$. How to build a high-quality $\beta^{\rm LB}$ efficiently 
is the topic of the next section.

\section{Natural-Norm SCM Lower Bound}
\label{sec:lower-bound}

The NNSCM  \cite{HKCHP}  constructs a decomposition of the (global) parameter domain 
\[\calD \equiv {\cup}_{k=1}^K \calD_{\mubar^k}\] by a greedy
approach. There is a ``control point'' $\mubar^k$ within each subdomain. 
Locally in each subdomain,  a linearized-in-parameter inf-sup formulation is utilized 
incorporating continuity information resulting in first order (in parameter) concavity. 
For the completeness of this paper and, in addition, due to that many ingredients of the NNSCM 
are adopted by our new cNNSCM, 
we detail the local and global aspects of this algorithm in the following two subsections respectively . 

\subsection{Local Approximation}
\label{subsec:local-approx}
The inf-sup numbers at the control points of these subdomains $\{\beta^\calN(\mubar^k)\}_{k=1}^K$ are 
calculated exactly and, at any other location, the ratio $\frac{\beta^\calN({\boldsymbol \mu})}{\beta^\calN(\mubar)}$ 
is approximated from below. Obviously the product of this lower bound and $\beta^\calN(\mubar^k)$ 
provides a lower bound for $\beta^\calN({\vmu})$.  Let us describe these two components briefly.

\subsubsection{From $\beta(\mubar)$ to $\beta({\boldsymbol \mu})$}

For a given subdomain $\calD_\mubar$ with control point $\mubar$, and $\forall  \, \vmu \in \calD_\mubar$, 
we define an inf-sup constant measured relative to a natural-norm \cite{SenNatNorm}:
$$\beta_{\mubar} ({\boldsymbol \mu})  = \inf_{w \in X} \sup_{v \in X} \frac{  a(w,v;{\boldsymbol \mu}) } { \|T^{\mubar} w\|_X \| v\|_X }; $$
and a lower bound for $\beta_\mubar({\boldsymbol \mu})$,
$$\bar{\beta}_\mubar({\boldsymbol \mu})  = \inf_{w \in X} \frac{ a(w,T^\mubar w;{\boldsymbol \mu}) } { \|T^\mubar w\|_X^2   } .$$ 
It can be easily shown that
$\bar{\beta}_\mubar({\boldsymbol \mu}) \le \beta_\mubar({\boldsymbol \mu})$ and that $\bar{\beta}_\mubar({\boldsymbol \mu})$
will be a good approximation
to $\beta_\mubar({\boldsymbol \mu})$ for $\vmu$ near $\mubar$ \cite{SenNatNorm}.
It is also straightforward to show that 
$\beta(\mubar) \bar{\beta}_{\mubar}({\boldsymbol \mu}) \le \beta({\boldsymbol \mu})$ which allows us to translate the lower bound
for $\bar{\beta}_{\mubar}({\boldsymbol \mu})$ into a lower  bound for $\beta({\boldsymbol \mu})$ given $\beta(\mubar)$.

\subsubsection{Reliable lower bound for $\frac{\beta({\boldsymbol \mu})}{\beta(\mubar)}$ through the SCM$^2$}

What remains of the local approximation is the application of the classical SCM$^2$ to construct lower and upper 
bounds for $\bar{\beta}_\mubar({\boldsymbol \mu})$. 
This is applicable by simply noting that
$$ \bar{\beta}_\mubar({\boldsymbol \mu}) = \inf_{y \in \calY_\mubar } \calJ (y;{\boldsymbol \mu}) , \mbox{  where  } \calJ (y;{\boldsymbol \mu}) = \sum_{q = 1}^Q \Theta_q({\boldsymbol \mu}) y_q $$
and
$$ \calY_\mubar  = \{ y \in {\mathbb R}^Q \; |\;  \exists\, w_y \in X \; {\rm s.t.}\;  y_q =
\frac{ a_q(w_y,T^\mubar w_y) } { \|T^\mubar w_y\|_X^2}, \quad 1 \le q \le Q \}.$$

However, for the completeness of this algorithm, let us provide the details of this procedure.
We first introduce the bounding box
$$ {\mathcal B}_\mubar = \prod_{q = 1}^Q \left[-\frac{\gamma_q}{\beta(\mubar)}, \frac{\gamma_q}{\beta(\mubar)} \right], \mbox{  where  } 
\gamma_q = \sup_{w \in X} \frac{ \| T_q w \|_X } {\| w \|_X }, \quad 1 \le q \le Q   ; $$
note that the $\gamma_q$ are independent of $\mubar$. 
Next, given the local SCM sample (whose construction, detailed in the next section, will be done in a greedy fashion)
$${\mathcal C}_\mubar =  \{ \hat{\mu}^1_\mubar, \ldots, \hat{\mu}^{J_\mubar}_\mubar \},$$ we can now define
$$ \calY^{\rm LB}_\mubar(\mu,\calC_{\mubar}) = \{ y \in \calB_\mubar \; | \; \sum_{q = 1}^Q \Theta_q(\mu') y_q \ge \bar{\beta}_\mubar(\mu'), \forall \mu' \in \calC_{\mubar}^{J_\mubar^{\rm nb}}(\mu) \},$$ 
and then the lower bound for $\bar{\beta}_\mubar ({\boldsymbol \mu})$ determined by $\calC_{\mubar}$ is obtained by solving the 
linear programming problem
\begin{equation}\label{SCMLB}
    \bar{\beta}^{\rm LB}_\mubar (\vmu;\calC_{\mubar}) = \inf_{y \in \calY^{\rm LB}_\mubar(\mu, \calC_{\mubar}) } \calJ (y;{\boldsymbol \mu}), \forall \vmu \in
\calD_\mubar .
\end{equation} 
Here, $ \calC_{\mubar}^{J_\mubar^{\rm nb}}(\mu)$ denotes the set of 
$J_\mubar^{\rm nb}$ points that are closest to $\mu$ within the set $\calC_{\mubar}$. 
To develop the upper bound, we simply introduce the set 
$$ \calY^{\rm UB}_\mubar(\calC_{\mubar}) = \{y^*_\mubar(\hat{\mu}^j_\mubar), 1 \leq j \leq J_{\mubar} \} \mbox{ where } y^*_\mubar({\boldsymbol \mu}) = \arg \min_{y\in\calY_\mubar}
\calJ(y;{\boldsymbol \mu});$$ 
and then define
\begin{equation}\label{SCMUB}
    \bar{\beta}^{\rm UB}_\mubar (\vmu;\calC_{\mubar}) = \inf_{y \in \calY^{\rm UB}_\mubar (\calC_{\mubar})} \calJ (y;{\boldsymbol \mu}), \quad \forall \vmu \in
\calD_\mubar .
\end{equation}

We want to make two remarks at this point:
\begin{itemize}
\item We have ${\mathcal Y}^{\rm UB}_\mubar(\calC_{\mubar}) \subset {\mathcal Y}_\mubar \subset \calY^{\rm LB}_\mubar(\mu, \calC_{\mubar})$ and hence \cite{HKCHP}
$\bar{\beta}^{\rm UB}_\mubar(\mu;\calC_{\mubar}) \ge \bar{\beta}_\mubar({\boldsymbol \mu}) \ge \bar{\beta}^{\rm LB}_\mubar (\mu;\calC_{\mubar}).$
\item The lower bound will only be useful
if $\bar{\beta}_{\mubar}({\boldsymbol \mu}) > 0$ over $\calD_\mubar$ which is, in general, 
not a consequence of $\beta({\boldsymbol \mu}) > 0, \forall \vmu \in \calD_{\mubar}$. 
We must thus adaptively divide the global parameter domain $\calD$ into 
subdomains $\calD_\mubar$ to ensure positivity. This is the subject of the next subsection. 
\end{itemize}

\subsection{Global Approximation: Greedy Sampling Procedure}
\label{subsec:global-approx}

We construct our domain decomposition $\calD = \cup_{k=1}^K \calD_{\mubar^k}$ by a greedy
approach. 
We first extend our lower and upper bounds of \eqref{SCMLB} and \eqref{SCMUB} to all $\vmu \in \calD$:
for a given $\mubar \in \calD$ and a finite sample $\calE \subset \calD$ we define
$$ g^{\rm LB}_{\mubar}(\vmu;\calE) = \inf_{y \in \calY^{\rm LB}_\mubar(\mu, \calE) } \calJ (y;{\boldsymbol \mu}), \forall \vmu \in \calD,$$ 
and
$$g^{\rm UB}_{\mubar}(\vmu;\calE) = \inf_{y \in \calY^{\rm UB}_\mubar(\calE) } \calJ (y;{\boldsymbol \mu}), \forall \vmu \in \calD.$$
This allows us to introduce an ``SCM$_{\rm R}$ quality control'' indicator.
$$\epsilon_{\mubar}(\vmu;\calE) \equiv
\frac { g^{\rm UB}_{\mubar}(\vmu;\calE)-
g^{\rm LB}_{\mubar}(\vmu;\calE) }
{g^{\rm UB}_{\mubar}(\vmu;\calE)}, \quad \forall \vmu \in \calD.$$
Here, the ``R'' in ``SCM$_{\rm R}$'' indicates that it is to control the ratio between $\beta({\boldsymbol \mu})$ and $\beta(\mubar)$.
Finally, we require a very rich train sample $\Xi \in \calD$, an SCM tolerance $\epsilon_\betabar \in (0,1)$, 
and an inf-sup tolerance function $\varphi(\mu,\mubar) \ge 0$ which is usually set to zero.

We are now ready to define the greedy algorithm in Algorithm \ref{alg:greedy}. 
\begin{algorithm}[h!]
  \caption{Natural-Norm SCM Greedy Algorithm}\label{alg:greedy}
  \begin{algorithmic}
\STATE {\bf 1.} Set $S = \{\mubar^1\}$, $k=1$; here $\mubar^1$ is an arbitrary point in $\Xi$;
\medskip
  \STATE {\bf 2.} Initialize $\calC_{\mubar^k} = \{ \empty \}$, $J_{\mubar^k} = 0$, $\calR_{\mubar^k} = \{\empty\}$, $\calR^*_{\mubar^k} = \{\empty\}$, and $\epsilon_{\rm max} = +\infty$
  \WHILE {$\calR_{\mubar^k}^* \backslash \calR_{\mubar^k} \neq\{\empty \}$ {\bf or}
$\epsilon_{\rm max} > \epsilon_{\betabar}$}
\STATE {\bf 2.1.} Set $\hat{\mu}_{\mubar^k}^{J_{\mubar^k}+1}$ to be
$\mubar^{k}$ if $J_{\mubar^k}=0$, and $\arg \max_{\mu \in \Xi}\epsilon_{\mubar^k}(\mu;\calC_{\mubar^k})$ otherwise.
\STATE {\bf 2.2.} Set $\calC_{\mubar^{k}} = \calC_{\mubar^{k}} \cup \{ \hat{\mu}_{\mubar^{k}}^{J_{\mubar^k}+1} \}$, $\calR_{\mubar^k} = \calR^*_{\mubar^k}$ and 
$\calR^*_{\mubar^k} = \{\mu \in \Xi: g^{\rm LB}_{\mubar^k}(\mu;\calC_{\mubar^k})> \varphi(\mu,\mubar^k)\}$,
$J_{\mubar^k} \leftarrow J_{\mubar^k} + 1$, $\epsilon_{\rm max} = \epsilon_{\mubar^k}(\mu;\calC_{\mubar^k})$.
\IF{$|\calR^*_{\mubar^k}| = |\Xi|$}
\STATE Break;
\ENDIF
\ENDWHILE
\medskip
\STATE {\bf 3.} Update (prune) $\Xi \leftarrow \Xi \backslash \calR_{\mubar^k}$;
\IF{$\Xi = \{ \empty \}$,}
\STATE Set $K = k$ and {\bf terminate};
\ELSE
\STATE Set $\mubar^{k+1} =  \arg \min_{\mu \in \Xi}
g^{\rm LB}_{\mubar^k}(\mu;\calC_{\mubar^k}),$ $S = S\cup\{\mubar^{k+1}\}$, $k = k + 1$, 
and {\bf goto} 2;
\ENDIF
  \end{algorithmic}
\end{algorithm}
The ``output" from the greedy procedure is the set of points $S = \{ \mubar^1,\ldots,\mubar^{K} \}$
and associated SCM sample sets $\calC_{\mubar^k}, 1 \le k \le K$. 
Several remarks regarding this algorithm follow:
\begin{itemize}
\item The set of points $\calR$ play the role
of temporary subdomains during the greedy construction. 
Observe that we declare the current subdomain/approximation complete  
(and move to the next subdomain) {\em only} when the trial sample offers no improvement in the positivity coverage 
{\em and} the trial sample is not required to satisfy our $\epsilon_{\betabar}$ SCM quality criterion. 
\item The improvement for a particular subdomain and identification of a new subdomain are based on different criteria. For the former
$\epsilon_{\mubar^k}(\vmu;\calC_{\mubar^k})$ is very effective: the $\arg \max$ will {\em avoid} $\vmu$ for which the upper bound
is negative and hence likely to lie outside the domain of relevance of $T^{\mubar^k}\!\!$, yet {\em favor} $\vmu$ for which the
current approximation is poor and hence (likely) to lie at the extremes of the domain of relevance of $T^{\mubar^k}\!\!$
- thus promoting optimal coverage. In contrast,
for the latter $g^{\rm LB}_{\mubar^k}$ is very effective: the arg min will look for the most negative value of the
lower bound - thus leaving the domain of relevance of $T^{\mubar^k}\!\!$. 
\end{itemize}

\noindent Finally, our global lower bound for $\beta({\boldsymbol \mu})$ is defined to be the maximum of those translated 
from each subdomain:
\begin{equation}\label{SCMgLB}
{\beta}^{\rm LB}({\boldsymbol \mu}) =  \max_{ k \in \{1 \ldots K \} } \beta(\mubar^{k}) g^{\rm LB}_{\mubar^{k}}( \vmu; \calC_{\mubar^{k}}).
\end{equation}

\section{Certified NNSCM}
\label{sec:cnnscm}

Our primary interest is in the lower bound $\beta^{\rm LB}({\boldsymbol \mu})$ as it is required for rigor in our
reduced basis {\it a posteriori} error estimator. However, the upper bound serves an important role 
in making sure the lower bound is not too pessimistic.
We note that NNSCM \cite{HKCHP} can ensure reasonable accuracy by choosing
an appropriate $\varphi(\mu,\mubar)$ in Algorithm \ref{alg:greedy}. 
However, there are usually parameters in $\varphi(\mu,\mubar)$ that we need to 
tune and there is no mechanism to easily control the true error of the lower bound for ${\beta({\boldsymbol \mu})}$.
 
Here, we develop an upper bound that can be constructed together with the natural-norm SCM
lower bound at marginal offline cost. 
To do that, we recall that for SCM$^2$,
$$ \beta^\calN({\boldsymbol \mu}) = \inf_{w\in X^\calN} \sup_{v \in X^\calN} \frac{a^\calN(w,v;{\boldsymbol \mu})}{\|w\|_{X}\,\|v\|_{X}}  \mbox{ and }
T^\mu w \equiv \sum_{q=1}^Q \Theta^q({\boldsymbol \mu})T^q w $$
to realize
\[
(\beta^{\mathcal N}({\boldsymbol \mu}))^2 = \inf_{w \in X^{\mathcal
N}}\sum_{q'=1}^Q\sum_{q''=q'}^Q(2 -
\delta_{q'q''})\Theta^{q'}({\boldsymbol \mu})\Theta^{q''}({\boldsymbol \mu})\frac{(T^{q'}w,T^{q''}w)_{X^{\mathcal
N}}}{\lVert w \rVert^2_{X^{\mathcal N}}}.
\]
Here, $\delta_{q'q''}$ is the Kronecker delta. Next, we identify
\begin{equation*}
(2 - \delta_{q'q''})\Theta^{q'}({\boldsymbol \mu})\Theta^{q''}({\boldsymbol \mu}), 1 \le q' \le
q'' \le Q \longmapsto \hat{\Theta}^q({\boldsymbol \mu}), 1 \le q \le \hat{Q}
\equiv \frac{Q(Q+1)}{2},
\end{equation*}
\begin{equation*}
\frac{(T^{q'}w,T^{q''}w)_{X^{\mathcal N}} + (T^{q''}w,T^{q'}w)_{X^{\mathcal N}}}{2},
1 \le q' \le q'' \le Q \longmapsto \hat{a}^\calN_q(w,v), 1 \le q \le
\hat{Q},
\end{equation*}
and obtain
\begin{equation}
\label{eq:nonco2co}
(\beta^{\mathcal N}({\boldsymbol \mu}))^2 \equiv
\inf_{w \in X^{\mathcal N}}
\sum_{q=1}^{\hat{Q}}\hat{\Theta}^q({\boldsymbol \mu})\frac{\hat{a}^\calN_q(w,w)}{\lVert
w \rVert^2_{X^{\mathcal N}}}.
\end{equation}

Having this interpretation, we simply introduce the set 
$$ w^{\rm UB}_\mubar(\calC_{\mubar}) = \{w^*_\mubar(\hat{\mu}^j_\mubar), 1 \leq j \leq J_{\mubar} \}$$ 
where $w^*_\mubar(\hat{\mu}^j_\mubar)$ is such that if we define
$y_q = \frac{ a_q(w,T^\mubar w) } { \|T^\mubar w\|_X^2}$ for $w = w^*_\mubar(\hat{\mu}^j_\mubar)$, $1 \le q \le Q$, we have
$\{y_1, \dots, y_Q\} = \arg \min_{y\in\calY_\mubar}\calJ(y;{\boldsymbol \mu}).$
We are now ready to define an upper bound for $(\beta^\calN({\boldsymbol \mu}))^2$,
 \begin{equation}
 \label{eq:SCMUB_Local}
 {\beta}^{\rm UB}_{\rm SCM^2}(\vmu; \mubar ) = \sqrt{\inf_{w \in w^{\rm UB}_\mubar (\calC_{\mubar})} \sum_{q=1}^{\hat{Q}}\hat{\Theta}^q({\boldsymbol \mu})\frac{\hat{a}^\calN_q(w,w)}{\lVert
w \rVert^2_{X^{\mathcal N}}}}, \quad \forall \mu \in \calD_\mubar .
\end{equation}

\begin{algorithm}[h!]
  \caption{Certified Natural-Norm SCM Greedy Algorithm}\label{alg:c_greedy}
  \begin{algorithmic}
\STATE {\bf 1.} Set $S = \{\mubar^1\}$, $k=1$; here $\mubar^1$ is an arbitrary point in $\Xi$;
\medskip
  \STATE {\bf 2.} Initialize $\calC_{\mubar^k} = \{ \empty \}$, $J_{\mubar^k} = 0$, $\calR_{\mubar^k} = \{\empty\}$, $\calR^*_{\mubar^k} = \{\empty\}$, and $\epsilon_{\rm max} = +\infty$
  \WHILE {$\calR_{\mubar^k}^* \backslash \calR_{\mubar^k} \neq\{\empty \}$ {\bf or}
$\epsilon_{\rm max} > \epsilon_{\betabar}$}
\STATE {\bf 2.1.} Set $\hat{\mu}_{\mubar^k}^{J_{\mubar^k}+1}$ to be
$\mubar^{k}$ if $J_{\mubar^k}=0$, and $\arg \max_{\mu \in \Xi}\epsilon_{\mubar^k}(\mu;\calC_{\mubar^k})$ otherwise.
\STATE {\bf 2.2.} Set $\calC_{\mubar^{k}} = \calC_{\mubar^{k}} \cup \{ \hat{\mu}_{\mubar^{k}}^{J_{\mubar^k}+1} \}$, $\calR_{\mubar^k} = \calR^*_{\mubar^k}$ and 
$\calR^*_{\mubar^k} = \{\mu \in \Xi: g^{\rm LB}_{\mubar^k}(\mu;\calC_{\mubar^k}) > 0\}$,
$J_{\mubar^k} \leftarrow J_{\mubar^k} + 1$, $\epsilon_{\rm max} = \epsilon_{\mubar^k}(\mu;\calC_{\mubar^k})$.
\IF{$|\calR^*_{\mubar^k}| = |\Xi|$}
\STATE Break;
\ENDIF
\ENDWHILE
\medskip
\STATE {\bf 3.} Update (prune) $\Xi \leftarrow \Xi \backslash \calR_{\mubar^k}$;
\IF{$\Xi = \{ \empty \}$,}
\IF{$\epsilon(\mu,\cdot) \le \epsilon_g$}
\STATE Set $K = k$ and {\bf terminate};
\ELSE
\STATE Set $k=1$, $S = \{\mubar^1\}$ with $\mubar^1 = \arg \max_{\mu \in \Xi} \epsilon({\boldsymbol \mu})$ and {\bf goto} 2.
\ENDIF
\ELSE
\STATE Set $\mubar^{k+1} =  \arg \min_{\mu \in \Xi}
g^{\rm LB}_{\mubar^k}(\mu;\calC_{\mubar^k}),$ $S = S\cup\{\mubar^{k+1}\}$, $k = k + 1$, 
and {\bf goto} 2;
\ENDIF
  \end{algorithmic}
\end{algorithm}

The global upper bound for $\beta({\boldsymbol \mu})$ is defined to be the minimum of 
${\beta}^{\rm UB}_{\rm SCM^2}(\vmu; \mubar )$ for different control points $\mubar$:
\begin{equation}\label{eq:SCMgUB}
{\beta}^{\rm UB}({\boldsymbol \mu}) =  \min_{ k \in \{1 \ldots K \} }  {\beta}^{\rm UB}_{\rm SCM^2}(\mu; \mubar_k ),
\end{equation}
and the ``SCM$_\beta$ quality control'' of the global lower bound
$$\epsilon({\boldsymbol \mu}) = \frac{{\beta}^{\rm UB}({\boldsymbol \mu})  - {\beta}^{\rm LB}({\boldsymbol \mu}) }{{\beta}^{\rm UB}({\boldsymbol \mu}) }.$$
We are now ready to state the certified NNSCM, Algorithm \ref{alg:c_greedy}. Here we introduce an additional 
tolerance $\epsilon_g$ which is to bound $\epsilon({\boldsymbol \mu})$ so that we have
$$(1 - \epsilon_g) {\beta}^{\rm UB}({\boldsymbol \mu}) < {\beta}^{\rm LB}({\boldsymbol \mu}) < {\beta}^\calN({\boldsymbol \mu}) < {\beta}^{\rm UB}({\boldsymbol \mu}).$$

This algorithm is very similar to Algorithm \ref{alg:greedy}. In addition to defining the global upper bound, it 
incorporates the mechanism of allowing for multiple rounds of domain decomposition which is not possible 
for NNSCM. In the context of the cNNSCM, NNSCM stops after the first round when the whole 
parameter domain is covered and, for each subdomain, the quality of the lower bound for the ratio 
$\frac{\beta(\vmu)}{\beta(\mubar)}$ has achieved the desired tolerance. 
On the other hand, the cNNSCM detects this, through monitoring the quality of $\beta(\vmu)$ measured by 
${\frac{\beta^{\rm UB}(\vmu) - \beta^{\rm LB}(\vmu)}{\beta^{\rm UB}(\vmu)}}$, and continue 
with more rounds of domain decomposition. For each $\mu \in \calD$, it is covered by one subdomain at each round 
making it possible to sharpen the lower bound $\beta^{\rm LB}(\vmu)$ in approximating $\beta(\vmu)$. 
Moreover, to start a new round, the size of $\epsilon({\boldsymbol \mu})$ is a good indicator 
for the need of a control point. Thus we set $\mubar^1 = \arg \max_{\mu \in \Xi} \epsilon({\boldsymbol \mu})$ 
to set the stage for the next domain decomposition.

Another important remark is that the increase in computational cost due to the $Q^2-$term expansion in 
\eqref{eq:nonco2co} is negligible. There is only one $\hat{Q}-$operation (a $\hat{Q}$ term summation in 
\eqref{eq:SCMUB_Local}) every time there is a control point identified or there is a new ${\rm SCM}_{\rm R}$ 
sample point added within a subdomain. This is negligible in comparison to the $\hat{Q}-$dependent cost in 
${\rm SCM}^2$ whose elimination is one critical contribution of NNSCM.

\section{Numerical Results}

\label{sec:numerical}

We test our implementation of the NNSCM and cNNSCM on the following two test problems:
\begin{subequations}\label{eq:example-pdes-2d}
\begin{align}\label{eq:awave2d}
  \left\{
  \begin{aligned}
    -u_{xx} - \mu^1 u_{yy} - \mu^2 u &= f(\x), & \x &\in \Omega \\
    u &= g, & \x &\in \partial \Omega
  \end{aligned}
  \right.
  & \hskip 10pt
  \mu \in \calD = [0.1,4] \times [0,2],
\end{align}
\begin{align}\label{eq:diff2d}
  \left\{
  \begin{aligned}
    (1+\mu^1 x) u_{xx} + (1+\mu^2 y) u_{yy} &= f(\x), & \x &\in \Omega \\
    u &= g, & \x &\in \partial \Omega
  \end{aligned}
  \right.
  & \hskip 10pt
  \mu \in \calD = [-0.99,0.99]^2.
\end{align}
\end{subequations}

\begin{figure}[htbp]
\begin{center}
\includegraphics[width=0.98\textwidth]{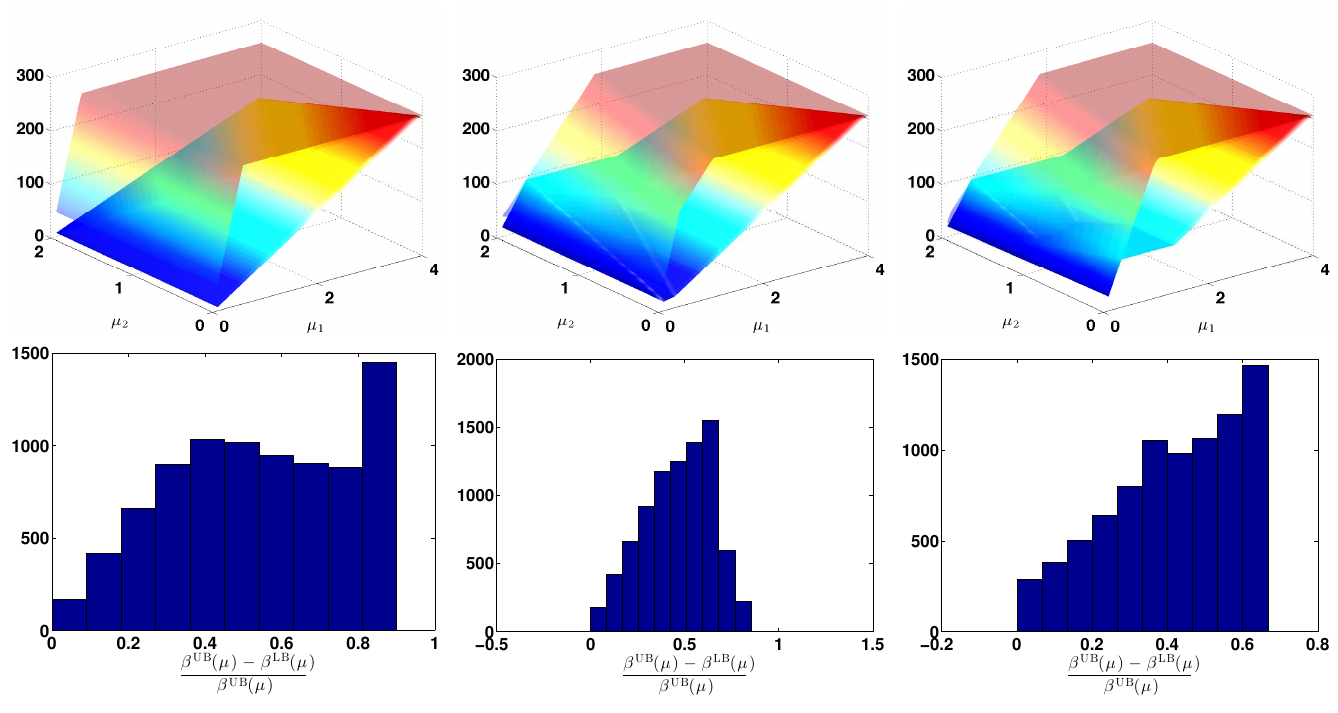}
\caption{Certified NNSCM for the first test problem. 
Plotted on the first row are the lower bounds and upper bounds (as functions of the parameter) after the $1^{\rm st}$, $2^{\rm nd}$, and $3^{\rm rd}$ (final) round of domain decomposition.  
On the bottom row is the histogram of ${\frac{\beta^{\rm UB}(\vmu) - \beta^{\rm LB}(\vmu)}{\beta^{\rm UB}(\vmu)}}$. 
Their maximums are $0.9, 0.85, 0.67$ from left to right. 
Notice that the original NNSCM would have stopped after the first column.}
\label{fig:P1result}
\end{center}
\end{figure}

The result for the first problem is shown in Figure \ref{fig:P1result}. We discretize the parameter domain 
by a $129 \times 65$ uniform grid, and the differential operator by the Pseudospectral collocation method 
\cite{TrefethenSpecBook}. 
We set $J_\mubar^{\rm nb} = 8$, $\epsilon_\betabar = \epsilon_g = 0.8$ 
(the later applicable to cNNSCM only). 
Plotted on the first row are $\beta^{\rm UB}(\vmu)$ overlaid 
to $\beta^{\rm LB}(\vmu)$. For the NNSCM, since the parameter domain $\calD$ is completely decomposed and the 
$\epsilon_\betabar$-condition is met on each subdomain, it will stop after the first column. As a result, the quality of 
the lower bound, measured by $\displaystyle{\max_\vmu \frac{\beta^{\rm UB}(\vmu) - \beta^{\rm LB}(\vmu)}{\beta^{\rm UB}(\vmu)}}$, 
is bad. On the other hand, the cNNSCM detects this 
and continue with two more rounds of domain decompositions gradually improving the quality of the lower bound. 
This is clearly visible on the graph and also shown by the decreasing of 
 $\displaystyle{\max_\vmu \frac{\beta^{\rm UB}(\vmu) - \beta^{\rm LB}(\vmu)}{\beta^{\rm UB}(\vmu)}}$. 
 To take a closer look at the quality of the lower bound, we plot the histogram of 
 $\displaystyle{\frac{\beta^{\rm UB}(\vmu) - \beta^{\rm LB}(\vmu)}{\beta^{\rm UB}(\vmu)}}$ 
on the second row. 
 It shows that the gap between the lower bound and upper bound is below the prescribed tolerance and the 
 quality of the lower bound is uniformly better than that by NNSCM.
 
\begin{figure}[htbp]
\begin{center}
\includegraphics[width=0.98\textwidth]{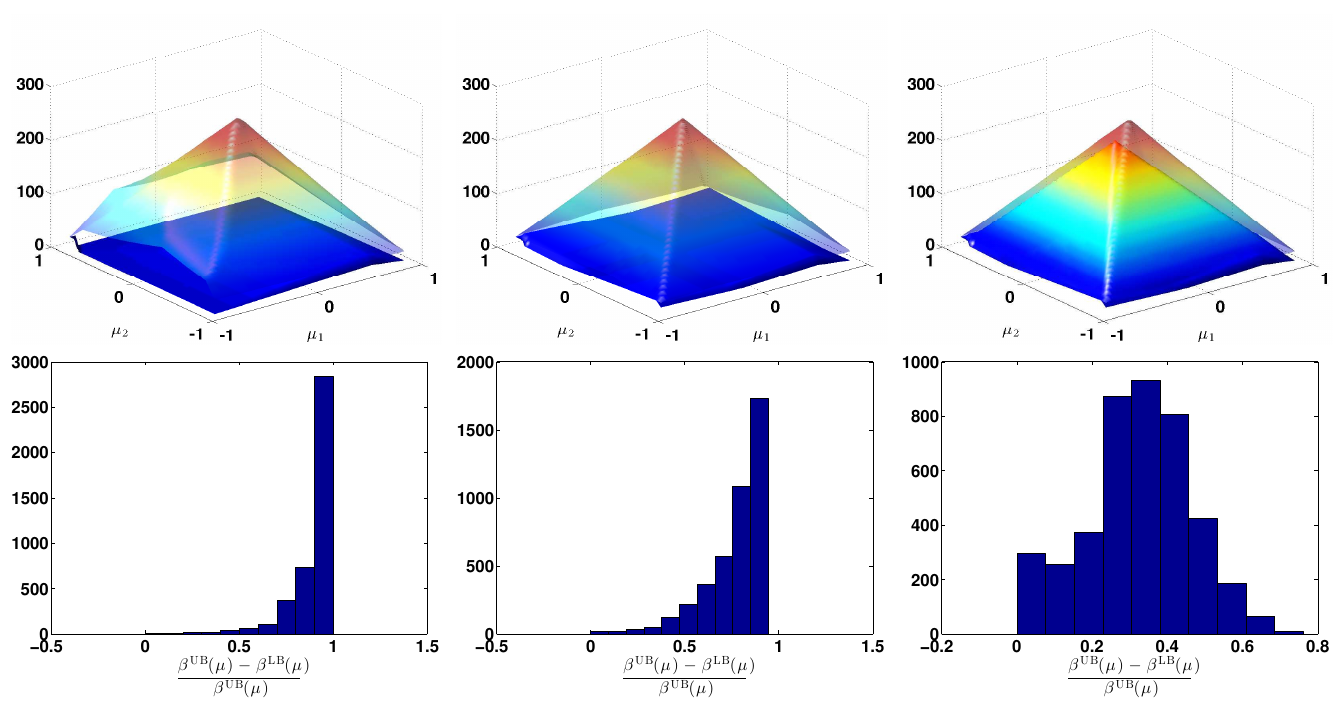}
\caption{Certified NNSCM for the second test problem.  
Plotted on the first row are the lower bounds and upper bounds (as functions of the parameter) after the $1^{\rm st}$, $4^{\rm th}$, and $6^{\rm th}$ (final) round of domain decomposition.  
On the bottom row is the histogram of ${\frac{\beta^{\rm UB}(\vmu) - \beta^{\rm LB}(\vmu)}{\beta^{\rm UB}(\vmu)}}$. 
Their maximums are $0.99986, 0.943, 0.761$ from left to right. 
Notice that the original NNSCM would have stopped after the first column.}
\label{fig:P2result}
\end{center}
\end{figure}

 The result for the second problem is shown in Figure \ref{fig:P2result}. The setting is the same other than 
 that the parameter domain is discretized by a $65 \times 65$ grid. This is a more challenging problem in the sense that 
 it becomes close to being degenerate at the four corners of the parameter domain. The poor quality of the NNSCM lower bound 
 is clearly shown by the picture and that $\displaystyle{\max_\vmu \frac{\beta^{\rm UB}(\vmu) - \beta^{\rm LB}(\vmu)}{\beta^{\rm UB}(\vmu)}} = 0.99986$ at convergence for NNSCM. 
 The cNNSCM has, again, improved it with a few more rounds of decompositions resulting in a lower bound that is 
 very close to the upper bound.


\section{Concluding remarks}
\label{sec:conclusion}

A rigorous and controllably tight lower bound for the stability parameter that is efficiently achievable 
is a critical part in the development of certified reduced
basis methods for parametrized partial differential equations. 
The available methodologies either 
suffer from significant computational cost or inferior tightness of the bound. 

In this paper, we have improved a recent novel approach combining two previous techniques by adding a mechanism 
through which the practitioners can control the tightness of the lower bound. It is achieved by building simultaneously an upper 
bound and shrinking the gap between the two through multiple domain decompositions. 
Numerical experiments demonstrate the effectiveness of the new approach and highlight its significant 
improvement over the Natural-Norm SCM.

\section*{Acknowledgements} The author would like to thank
Professor Hesthaven, Jan from EPFL for encouragement and helpful discussions during the process of this project.


\begin{thebibliography}{10}

\bibitem{Adams1975Book}
R.A. Adams.
\newblock {\em Sobolev Spaces}.
\newblock Pure and applied mathematics. Academic Press, 1975.

\bibitem{Almroth_Stern_Brogan}
B.~O. Almroth, P.~Stern, and F.~A. Brogan.
\newblock Automatic choice of global shape functions in structural analysis.
\newblock {\em AIAA Journal}, 16:525--528, May 1978.

\bibitem{Barrault_Nguyen_Maday_Patera}
M.~Barrault, N.~C. Nguyen, Y.~Maday, and A.~T. Patera.
\newblock An ``empirical interpolation'' method: Application to efficient
  reduced-basis discretization of partial differential equations.
\newblock {\em C. R. Acad. Sci. Paris, S\'erie I}, 339:667--672, 2004.

\bibitem{CHMR-Cras}
Y.~Chen, J.~S. Hesthaven, Y.~Maday, and J.~Rodr\'{i}guez.
\newblock A monotonic evaluation of lower bounds for inf-sup stability
  constants in the frame of reduced basis approximations.
\newblock {\em C. R. Acad. Sci. Paris, Ser. I}, 346:1295--1300, 2008.

\bibitem{CHMR-M2an}
Y.~Chen, J.~S. Hesthaven, Y.~Maday, and J.~Rodr\'{i}guez.
\newblock Improved successive constraint method based a posteriori error
  estimate for reduced basis approximation of 2d maxwell's problem.
\newblock {\em M2AN}, 43:1099--1116, 2009.

\bibitem{Deparis08}
S.~Deparis.
\newblock Reduced basis error bound computation of parameter-dependent
  navier-stokes equati the natural norm approach.
\newblock {\em SIAM J. Numer. Anal.}, 46(4):2039--2067, 2008.

\bibitem{Fink_Rheinboldt_1}
J.~P. Fink and W.~C. Rheinboldt.
\newblock On the error behavior of the reduced basis technique for nonlinear
  finite element approximations.
\newblock {\em Z. Angew. Math. Mech.}, 63(1):21--28, 1983.

\bibitem{Grepl_Maday_Nguyen_Patera}
M.~A. Grepl, Y.~Maday, N.~C. Nguyen, and A.~T. Patera.
\newblock Efficient reduced-basis treatment of nonaffine and nonlinear partial
  differential equations.
\newblock {\em Mathematical Modelling and Numerical Analysis}, 41(3):575--605,
  2007.

\bibitem{HaasdonkOhlberger}
B. Haasdonk, M. Ohlberger. \newblock Reduced basis method for finite volume approximations
of parametrized linear evolution equations. \newblock {\em M2AN Math. Model. Numer. Anal.} 42: 277--302, 2008.

\bibitem{HKCHP}
D.B.P. Huynh, D.J. Knezevic, Y.~Chen, J.S. Hesthaven, and A.T. Patera.
\newblock A natural-norm successive constraint method for inf-sup lower bounds.
\newblock {\em CMAME}, 199:1963--1975, 2010.

\bibitem{HuynhSCM}
D.B.P. Huynh, G.~Rozza, S.~Sen, and A.T. Patera.
\newblock A successive constraint linear optimization method for lower bounds
  of parametric coercivity and inf-sup stability constants.
\newblock {\em C. R. Acad. Sci. Paris, S$\acute{e}$rie I.}, 345:473 -- 478,
  2007.

\bibitem{Machielis_Maday_Oliveira_Patera_Rovas}
L.~Machiels, Y.~Maday, I.~B. Oliveira, A.~T. Patera, and D.~V. Rovas.
\newblock Output bounds for reduced-basis approximations of symmetric positive
  definite eigenvalue problems.
\newblock {\em C. R. Acad. Sci. Paris S\'er. I Math.}, 331(2):153--158, 2000.

\bibitem{Maday}
Y.~Maday.
\newblock Reduced basis method for the rapid and reliable solution of partial
  differential equations.
\newblock In {\em International {C}ongress of {M}athematicians. {V}ol. {III}},
  pages 1255--1270. Eur. Math. Soc., Z\"urich, 2006.

\bibitem{Maday_Patera_Rovas}
Y.~Maday, A.~T. Patera, and D.~V. Rovas.
\newblock A blackbox reduced-basis output bound method for noncoercive linear
  problems.
\newblock In {\em Nonlinear partial differential equations and their
  applications. Coll\`ege de France Seminar, Vol. XIV (Paris, 1997/1998)},
  volume~31 of {\em Stud. Math. Appl.}, pages 533--569. North-Holland,
  Amsterdam, 2002.

\bibitem{NguyenVeroyPatera2005}
N.C. Nguyen, K.~Veroy, and A.~T. Patera.
\newblock Certified real-time solution of parametrized partial differential
  equations.
\newblock In Sidney Yip, editor, {\em Handbook of Materials Modeling}, pages
  1529--1564. Springer Netherlands, 2005.

\bibitem{Noor_Peters}
A.~K. Noor and J.~M. Peters.
\newblock Reduced basis technique for nonlinear analysis of structures.
\newblock {\em AIAA Journal}, 18(4):455--462, April 1980.

\bibitem{Porsching}
T.~A. Porsching.
\newblock Estimation of the error in the reduced basis method solution of
  nonlinear equations.
\newblock {\em Math. Comp.}, 45(172):487--496, 1985.

\bibitem{Prudhomme_Rovas_Veroy_Maday_Patera_Turinici}
C.~Prud'homme, D.~Rovas, K.~Veroy, Y.~Maday, A.~T. Patera, and G.~Turinici.
\newblock Reliable real-time solution of parametrized partial differential
  equations: Reduced-basis output bound methods.
\newblock {\em Journal of Fluids Engineering}, 124(1):70--80, March 2002.

\bibitem{QuarteroniVallie2008book}
A.~Quarteroni and A.~Valli.
\newblock {\em Numerical Approximation of Partial Differential Equations}.
\newblock Springer Series in Computational Mathematics. Springer, 2008.

\bibitem{Rozza_Huynh_Patera}
G.~Rozza, D.B.P. Huynh, and A.T. Patera.
\newblock Reduced basis approximation and a posteriori error estimation for
  affinely parametrized elliptic coercive partial differential equations:
  Application to transport and continuum mechanics.
\newblock {\em Arch Comput Methods Eng}, 15(3):229--275, 2008.

\bibitem{SenNatNorm}
S.~Sen, K.~Veroy, D.B.P. Huynh, S.~Deparis, N.C. Nguyen, and A.T. Patera.
\newblock ``{Natural norm}'' a posteriori error estimators for reduced basis
  approximations.
\newblock {\em J. Comput. Phys.}, 217(1):37 -- 62, 2006.

\bibitem{TrefethenSpecBook}
L.~N. Trefethen.
\newblock {\em Spectral methods in {MATLAB}}, volume~10 of {\em Software,
  Environments, and Tools}.
\newblock Society for Industrial and Applied Mathematics (SIAM), Philadelphia,
  PA, 2000.
  
\bibitem{Vallaghe}
S. Vallagh$\acute{\rm e}$, A. Le Hyaric, M. Fouquembergh, and C. Prud'homme. \newblock A successive constraint
method with minimal offline constraints for lower bounds of parametric coercivity constant.
\newblock Preprint: hal-00609212, hal.archives-ouvertes.fr

\bibitem{Zhang2011}  
S. Zhang. \newblock Efficient greedy algorithms for successive constraints methods with high-dimensional
parameters. \newblock {\em Brown Division of Applied Math Scientific Computing}, Tech Report, 23, 2011.





\end{thebibliography}
\end{document}